\documentclass[preprint,12pt]{article}
\usepackage{amsthm}
\usepackage{amssymb}
\usepackage{amsmath}
\usepackage{hyperref}
\usepackage{array}
\usepackage{booktabs}
\usepackage{color}
\usepackage{graphicx}
\usepackage[T1]{fontenc}
\usepackage{standalone}
\usepackage{tikz}
\usetikzlibrary{calc}

\theoremstyle{definition}
\newtheorem{problem}{Problem}
 
\begin{document}

\title{Searching for regular, triangle-distinct graphs\thanks{%
       E-mails: dragan.stevanovic@aasu.edu.kw, mohammad.ghebleh@ku.edu.kw, gilles.caporossi@hec.ca,
                vambat@gmail.com, sanja.stevanovic@mi.sanu.ac.rs}}

\author{Dragan Stevanovi\'c\thanks{On leave from the Mathematical Institute of Serbian Academy of Sciences and Arts.}, \\
        Abdullah Al Salem University, Kuwait
\and
        Mohammad Ghebleh, \\
        Kuwait University, Kuwait
\and        
        Gilles Caporossi, \\
        HEC Montr\'eal, Canada
\and        
        Ambat Vijayakumar, \\
        Cochin University of Science and Technology, India
\and        
        Sanja Stevanovi\'c\thanks{SS is supported by the Science Fund of the Republic of Serbia, 
             grant \#6767, Lazy walk counts and spectral radius of threshold graphs---LZWK.}, \\
        Mathematical Institute SANU, Serbia}

\maketitle

\begin{abstract}
The triangle-degree of a vertex~$v$ of a simple graph~$G$ is the number of triangles in~$G$ that contain~$v$.
A simple graph is triangle-distinct if all its vertices have distinct triangle-degrees.
Berikkyzy et al.\ [Discrete Math. 347 (2024) 113695] recently asked
whether there exists a regular graph that is triangle-distinct.
Here we showcase the examples of regular, triangle-distinct graphs with orders between 21 and~27,
and report on the methodology used to find them.


\end{abstract}

Let $G=(V,E)$ be a simple graph with $n=|V|$ vertices.
The triangle-degree $t(v)$ of vertex~$v\in V$ is the number of triangles in~$G$ that contain~$v$.
While every simple graph must contain two vertices of equal degree
(since it cannot simultaneously contain vertices of degrees 0 and~$n-1$),
there exist {\em triangle-distinct} graphs in which all vertices have distinct triangle-degrees.
We will abbreviate triangle-distinct graphs as TD graphs in the sequel.
According to Berikkyzy et al.~\cite{Beri},
Erd\"os and Trotter initiated the study of TD graphs
by asking for bounds on the number of edges in such graphs.
While Nair and Vijayakumar studied properties of triangle-degrees
in graphs and their complements~\cite{Vijay1} and in graph compositions~\cite{Vijay2},
it seems that \cite{Beri} is actually the first published paper that specifically mentions TD graphs.
Besides constructing an infinite family of TD graphs and
providing a non-trivial lower bound on the number of edges in TD graphs,
Berikkyzy et al.~\cite{Beri} also posed the following problem.

\begin{problem}[\cite{Beri}]
Does there exist a regular graph that is triangle-distinct?
\end{problem}

It is not hard to find concrete examples of general TD graphs.
Berikkyzy et al.~\cite{Beri} noted that the smallest TD graph has seven vertices only,
and our enumeration of simple graphs with up to 11~vertices yields the following results:
\begin{center}
\begin{tabular}{rrrc}
\toprule
Order & Simple graphs & TD graphs & Percentage (\%) \\
\midrule
7  &         1,044 &         1 & 0.096 \\
8  &        12,346 &        31 & 0.251 \\
9  &       274,668 &       924 & 0.336 \\
10 &    12,005,168 &    49,088 & 0.409 \\
11 & 1,018,997,864 & 4,389,900 & 0.431 \\
\bottomrule
\end{tabular}
\end{center}
However, while TD graphs do appear to become more common with increasing order,
none of TD graphs with up to 11 vertices turned out to be regular.

There does not seem to be an obvious reason why a regular TD graph should not exist.
If $A$ denotes the adjacency matrix of an $r$-regular TD graph~$G$ with $n$~vertices,
then the diagonal entries of~$A^2$ represent the vertex degrees,
while the diagonal entries of~$A^3$ represent twice the triangle-degrees.
Hence regularity is ensured by requiring all diagonal entries of~$A^2$ to be equal to~$r$,
while triangle-distinctness is ensured by requiring all diagonal entries of~$A^3$ to be mutually distinct.
If we choose $n$ distinct numbers between 0 and~${r\choose 2}$ to represent triangle-degrees
(provided that ${r\choose 2}\geq n-1$, i.e., that $r>\sqrt{2n}$, see also \cite[Theorem~5]{Beri}),
then these requirements together yield a system of $2n$~polynomial (i.e., quadratic and cubic) equations
in terms of the ${n\choose 2}$ entries of the upper part of~$A$.
While the current knowledge of solving systems of polynomial equations is 
not as advanced as we might like it (see, e.g., \cite{Sturm}),
our intuition said that being free to choose the values of ${n\choose 2}$~binary variables
should be enough to satisfy those $2n$~polynomial equations.

Led by this intuition,
we next ran an exhaustive search among regular graphs with up to 15 vertices,
which were generated by Brendan McKay's program geng~\cite{nauty}.
The largest set here was that of 1,470,293,675 8-regular graphs on 15 vertices,
but still without an example of a regular TD graph.

Apparently, this means that regular TD graphs live on larger numbers of vertices,
which cannot (yet) be exhaustively searched in a reasonable amount of time.
In such situation we had to resort to optimisation approach,
where the object that we aim for is represented
as an element of a search space 
in which an appropriately defined objective function reaches the optimum value.
The search space is obviously the set $\mathcal{RG}_{n,r}$ of all $r$-regular graphs on $n$~vertices,
for some fixed values of $n\geq 16$ and $r<n$.
There are at least two straightforward choices for the objective function~$f(G$):
\begin{itemize}
\item $f_1(G)$ is the cardinality of the set $\{t(u)\colon u\in V\}$, and
      then TD graphs are those graphs~$G$ in~$\mathcal{RG}_{n,r}$ for which $f_1(G)$ reaches~$n$ as its maximum value;
\item $f_2(G)$ is the cardinality of the set $\{\{u,v\}\colon t(u)=t(v), u\neq v, u,v\in V\}$, and
      then TD graphs are those graphs~$G$ in~$\mathcal{RG}_{n,r}$ for which $f_2(G)$ reaches~$0$ as its minimum value.
\end{itemize}
Additionally, 
in order to be able to apply local search methods
we have to introduce a neighbourhood relation on the set~$\mathcal{RG}_{n,r}$.
Let $u,v,s$ and $t$ be four distinct vertices of the graph $G=(V,E)$ 
such that $us,vt\in E$, while $ut,vs\notin E$.
The {\em edge switching} consists of 
deleting the edges $us$ and $vt$ and adding the edges $ut$ and $vs$,
so that after edge switching the new graph becomes $G-\{us,vt\}+\{ut,vs\}$.
Note that edge switching does not change degrees of the involved vertices,
so that if $G\in\mathcal{RG}_{n,r}$ then $G-\{us,vt\}+\{ut,vs\}\in\mathcal{RG}_{n,r}$ as well.
Our neighbourhood relation on~$\mathcal{RG}_{n,r}$ is then described in terms of edge switchings:
two graphs $G,H\in\mathcal{RG}_{n,r}$ are neighbours
if one can be obtained from the other by some edge switching.

We first implemented the greedy search which, for given $n$ and~$r$,
starts with a random $r$-regular graph on $n$~vertices (see \cite[Section~1.4]{java}) 
and at each next step
it enumerates all feasible edge switchings in the current graph and 
then moves to the neighbouring graph that mostly increases the value of~$f_1$.
This, however, did not lead to promising results
because the objective function~$f_1$ has only $n$~distinct values,
so that the greedy search easily goes astray and finishes searching way too prematurely.

The second greedy search, which moves to the neighbouring graph that mostly decreases the value of~$f_2$,
was much more encouraging.
It occasionally led to examples of regular graphs having {\em only one} vertex pair with equal triangle-degrees,
but most often it got stuck at local minima having two or three such vertex pairs.
Experiments with this greedy search also pointed out 
the most promising choice for the degree~$r$ to be around~$n/2$.

Motivated by this ``almost there'' situation,
we moved on to the variable neighbourhood search (VNS)
to find the minimisers of~$f_2$ in~$\mathcal{RG}_{n,r}$.
VNS, proposed by Mladenovi\'c and Hansen~\cite{vns}, is a clever buildup over the greedy search.
It first runs the greedy search to locate one local minimum~$G^\ast$ for~$f_2$ in~$\mathcal{RG}_{n,r}$.
In each subsequent iteration,
VNS ``shakes'' the current local minimum~$G^\ast$
by applying a sequence of $k$~random edge switchings to it to obtain a new graph~$H$,
and then starts the greedy search from~$H$.
If this greedy search ends up with the same local minimum value~$f_2(G^\ast)$,
VNS increases the value of~$k$ by one and repeats the shaking with the larger value of~$k$.
However, if the greedy search manages to find a better value for~$f_2$,
this smaller local minimum is set as the new value of~$G^\ast$ and $k$ is reset to one,
so that the future shakings start to explore bigger and bigger neighbourhoods of the new local minimum.
The hypothesis of VNS is that with enough many random shakings
it will be able to climb out of the valley imposed by any current local minimum and
that it will eventually reach the global minimum.
Note that VNS is essentially an endless procedure---%
unless you know how to recognise the global minimum you are looking for,
the only way to stop it is to put limits
either on the maximum shaking number~$k$ or on the time spent applying VNS.

VNS is the main ingredient of the program AutoGraphiX~\cite{agx1}, now in its third generation~\cite{agx3},
which has already served to verify or disprove many conjectures in graph theory~\cite{agxsurv1,agxsurv2}.
However, AutoGraphiX cannot be restricted to consider regular graphs only,
so we separately implemented VNS for $\mathcal{RG}_{n,r}$,
which is freely available at \url{https://github.com/dragance106/VNSregular}
and could be used as a base for other computational studies on regular graphs.

However, it turned out that even VNS could not find regular graphs with $f_2(G)<1$.
This called for the introduction of yet another objective function~$f_3$
that has even more diverse values over the set~$\mathcal{RG}_{n,r}$.
What we needed is a function that severely penalises the appearance of vertices with equal triangle-degrees,
so that a TD graph is easily recognisable from a small value of~$f_3$.
Supposing that the triangle-degrees of vertices in~$G$ are
ordered as $t_1\geq t_2\geq\cdots\geq t_n$,
we opted for the objective function
$$
f_3(G) = \sum_{i=1}^{n-1} \frac1{t_i-t_{i+1}+\frac1n}.
$$
In this way,
whenever there is a pair of vertices with equal triangle-degrees,
the value of~$f_3$ will be larger than~$n$,
while in TD graphs
all summands above will be smaller than one,
so that the value of~$f_3$ will be smaller than~$n$.

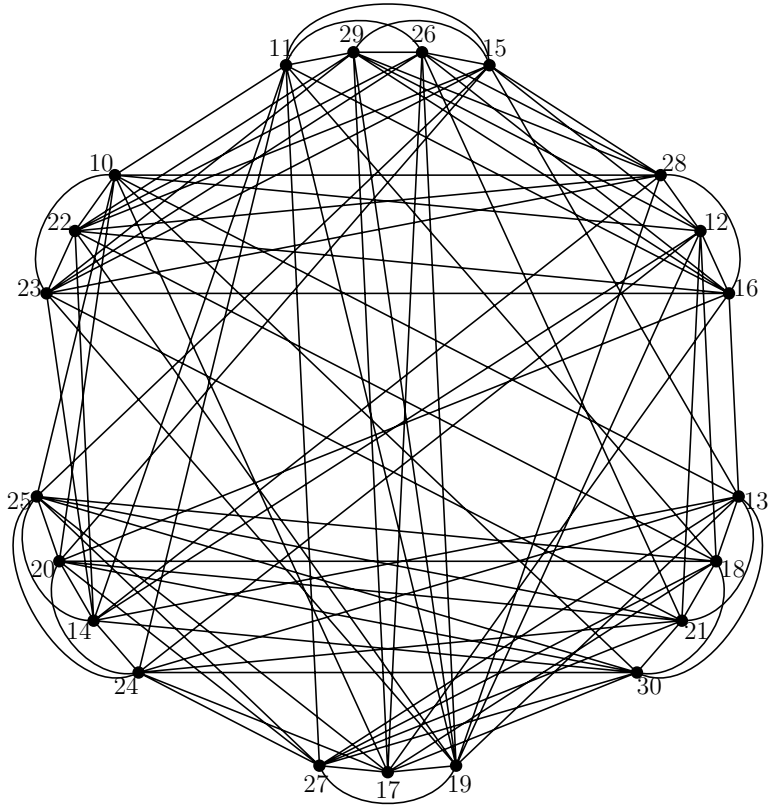
\begin{figure}[ht]
\begin{center}
\def\vtx{.08}
\def\scl{.60}
\def\nscl{.90}
\def\u{8}

\begin{tikzpicture}[rotate=180,scale=\scl,yscale=-1,semithick,inner sep=0pt]
\node (v0) at (0,0) {};
\node (v1) at (259.090909:\u) {};
\node (v2) at (270.000000:\u) {};
\node (v3) at (280.909091:\u) {};
\node (v4) at (313.636364:\u) {};
\node (v5) at (324.545454:\u) {};
\node (v6) at (335.454546:\u) {};
\node (v7) at (346.363636:\u) {};
\node (v8) at (379.090909:\u) {};
\node (v9) at (390.000000:\u) {};
\node (v10) at (400.909091:\u) {};
\node (v11) at (433.636364:\u) {};
\node (v12) at (444.545454:\u) {};
\node (v13) at (455.454546:\u) {};
\node (v14) at (466.363636:\u) {};
\node (v15) at (499.090909:\u) {};
\node (v16) at (510.000000:\u) {};
\node (v17) at (520.909091:\u) {};
\node (v18) at (553.636364:\u) {};
\node (v19) at (564.545454:\u) {};
\node (v20) at (575.454546:\u) {};
\node (v21) at (586.363636:\u) {};

\draw (v1)--(v8) (v1)--(v9) (v1)--(v11) (v1)--(v12) (v1)--(v13) (v1)--(v15) (v1)--(v16) (v1)--(v18) (v2)--(v4) (v2)--(v7) (v2)--(v10) (v2)--(v12) (v2)--(v13) (v2)--(v17) (v2)--(v19) (v2)--(v21) (v3)--(v4) (v3)--(v6) (v3)--(v7) (v3)--(v11) (v3)--(v18) (v3)--(v19) (v3)--(v20) (v3)--(v21) (v4)--(v11) (v4)--(v16) (v4)--(v18) (v4)--(v20) (v4)--(v21) (v5)--(v8) (v5)--(v9) (v5)--(v11) (v5)--(v15) (v5)--(v16) (v5)--(v18) (v5)--(v21) (v6)--(v10) (v6)--(v14) (v6)--(v17) (v6)--(v19) (v6)--(v20) (v6)--(v21) (v7)--(v10) (v7)--(v14) (v7)--(v19) (v7)--(v20) (v7)--(v21) (v8)--(v12) (v8)--(v13) (v8)--(v14) (v8)--(v15) (v8)--(v17) (v8)--(v20) (v9)--(v12) (v9)--(v13) (v9)--(v14) (v9)--(v15) (v9)--(v17) (v9)--(v19) (v10)--(v11) (v10)--(v15) (v10)--(v16) (v10)--(v18) (v10)--(v21) (v11)--(v17) (v11)--(v19) (v12)--(v15) (v12)--(v16) (v12)--(v17) (v13)--(v15) (v13)--(v17) (v13)--(v20) (v14)--(v15) (v14)--(v16) (v14)--(v18) (v16)--(v19) (v16)--(v20) (v17)--(v18);

\draw (v1)--(v2)--(v3);
\draw (v1) to[bend right=65] (v3);

\draw (v4)--(v5)--(v6)--(v7);
\draw (v4) to[bend right=65] (v6);
\draw (v5) to[bend right=65] (v7);
\draw (v4) to[bend right=85] (v7);

\draw (v8)--(v9)--(v10);
\draw (v8) to[bend right=65] (v10);

\draw (v11)--(v12)--(v13)--(v14);
\draw (v11) to[bend right=65] (v13);
\draw (v12) to[bend right=65] (v14);
\draw (v11) to[bend right=85] (v14);

\draw (v15)--(v16)--(v17);
\draw (v15) to[bend right=65] (v17);

\draw (v18)--(v19)--(v20)--(v21);
\draw (v18) to[bend right=65] (v20);
\draw (v19) to[bend right=65] (v21);
\draw (v18) to[bend right=85] (v21);

\foreach \i in {1,...,21}{
  \draw[fill] (v\i) circle (1.5*\vtx);
}

\def\labelscale{.8}
\node[scale=\labelscale] at ($1.05*(v1)$) {$19$};
\node[scale=\labelscale] at ($1.05*(v2)$) {$17$};
\node[scale=\labelscale] at ($1.05*(v3)$) {$27$};
\node[scale=\labelscale] at ($1.05*(v4)$) {$24$};
\node[scale=\labelscale] at ($1.05*(v5)$) {$14$};
\node[scale=\labelscale] at ($1.05*(v6)$) {$20$};
\node[scale=\labelscale] at ($1.05*(v7)$) {$25$};
\node[scale=\labelscale] at ($1.05*(v8)$) {$23$};
\node[scale=\labelscale] at ($1.05*(v9)$) {$22$};
\node[scale=\labelscale] at ($1.05*(v10)$) {$10$};
\node[scale=\labelscale] at ($1.05*(v11)$) {$11$};
\node[scale=\labelscale] at ($1.05*(v12)$) {$29$};
\node[scale=\labelscale] at ($1.05*(v13)$) {$26$};
\node[scale=\labelscale] at ($1.05*(v14)$) {$15$};
\node[scale=\labelscale] at ($1.05*(v15)$) {$28$};
\node[scale=\labelscale] at ($1.05*(v16)$) {$12$};
\node[scale=\labelscale] at ($1.05*(v17)$) {$16$};
\node[scale=\labelscale] at ($1.05*(v18)$) {$13$};
\node[scale=\labelscale] at ($1.05*(v19)$) {$18$};
\node[scale=\labelscale] at ($1.05*(v20)$) {$21$};
\node[scale=\labelscale] at ($1.05*(v21)$) {$30$};

\end{tikzpicture}
\end{center}
\caption{An example of a 10-regular, triangle-distinct graph on 21 vertices.}
\label{fig-1}
\end{figure}

Using VNS to minimise~$f_3$ over $\mathcal{RG}_{n,r}$,
we were finally able to identify 11 examples of regular TD graphs
with the values of~$n$ ranging from 21 to 27 and $r$ ranging between 10 and~13.
This has taken several days on two everyday laptops
running several instances of VNS search on separate processor cores in parallel.
Fig.~\ref{fig-1} shows one example of a TD graph in $\mathcal{RG}_{21,10}$.

Adjacency lists of this and the remaining examples that we found are provided in~\ref{appx}.
Note that the complement of a regular TD graph is again a TD graph,
due to \cite[Corollary~2.5]{Vijay1}
which states that for a vertex~$u$ of an $r$-regular graph~$G$ of order~$n$
$$
t_{\overline G}(u) = {n-1\choose 2} - \frac32r(n-r-1) - t_G(u),
$$
so if all the values~$t_G(u)$ are mutually distinct,
then so are all the values~$t_{\overline G}(u)$.
The examples we found with VNS are neither self-complementary nor form complementary pairs,
so that each example listed in~\ref{appx} provides its complement as an additional example of a regular TD graph.

To conclude,
we have provided here a positive answer to the question of Berikkyzy et al.~\cite{Beri}
whether there exist regular, triangle-distinct graphs.
They do, but they happen to have larger number of vertices and to be relatively rare,
which makes them largely inaccessible to exhaustive search methods.
Optimisation methods are capable of finding examples of such graphs,
but even then caution is needed in setting up the proper search objective,
as illustrated here through our experience with the ill-chosen initial objective functions $f_1$ and~$f_2$.
We leave to the interested reader the task of trying to spot
some structure in the examples of regular TD graphs given in the appendix.
Note, however, that this is not necessarily a straightforward task 
because TD graphs are necessarily asymmetric
(any nontrivial automorphism implies equal triangle-degrees for any two similar vertices),
while it is only natural for us to first look for the presence of symmetry.

\section*{Acknowledgement}

DS is on leave from the Mathematical Institute of the Serbian Academy of Sciences and Arts.
Research of SS was supported by the Science Fund of the Republic of Serbia,
grant \#6767, Lazy walk counts and spectral radius of threshold graphs---LZWK.

\section*{Declaration of competing interest}

The authors declare that there are no conflicts of interest regarding this paper.

\appendix

\newpage
\section{Adjacency lists of regular, triangle-distinct graphs}
\label{appx}

Here we give the adjacency lists of 11 regular TD graphs,
including the one shown in Fig.~\ref{fig-1},
which we found by minimising the objective function~$f_3$ with~VNS.
Each graph is presented in a separate table, where the $i$th row presents
the list of neighbours of the vertex $i$, as well as its triangle-degree.

\subsection*{\rm Regular triangle-distinct graph \#$1$: $n=21$, $r=10$ (See Fig.~\ref{fig-1})}

\[\begin{array}{c<{\hspace{1.5em}}cccccccccc>{\hspace{1em}}c}
\toprule
v & \multicolumn{10}{c}{\text{Neighbours of }v} & t(v)\\
\midrule
\mathbf{ 1} &  2 &  3 &  4 &  8 & 11 & 13 & 14 & 16 & 19 & 21 & 10\\
\mathbf{ 2} &  1 &  5 &  6 &  7 &  9 & 10 & 15 & 17 & 18 & 20 & 11\\
\mathbf{ 3} &  1 &  5 &  6 &  7 &  9 & 10 & 12 & 15 & 19 & 20 & 12\\
\mathbf{ 4} &  1 &  5 &  6 &  7 &  9 & 10 & 12 & 15 & 18 & 21 & 13\\
\mathbf{ 5} &  2 &  3 &  4 & 11 & 13 & 14 & 15 & 16 & 19 & 21 & 14\\
\mathbf{ 6} &  2 &  3 &  4 & 11 & 13 & 14 & 16 & 17 & 19 & 20 & 15\\
\mathbf{ 7} &  2 &  3 &  4 &  8 & 11 & 13 & 14 & 17 & 19 & 20 & 16\\
\mathbf{ 8} &  1 &  7 &  9 & 10 & 15 & 16 & 17 & 18 & 20 & 21 & 17\\
\mathbf{ 9} &  2 &  3 &  4 &  8 & 11 & 12 & 13 & 16 & 18 & 21 & 18\\
\mathbf{10} &  2 &  3 &  4 &  8 & 13 & 14 & 17 & 18 & 19 & 20 & 19\\
\mathbf{11} &  1 &  5 &  6 &  7 &  9 & 12 & 15 & 16 & 18 & 21 & 20\\
\mathbf{12} &  3 &  4 &  9 & 11 & 14 & 15 & 16 & 17 & 18 & 21 & 21\\
\mathbf{13} &  1 &  5 &  6 &  7 &  9 & 10 & 14 & 17 & 19 & 20 & 22\\
\mathbf{14} &  1 &  5 &  6 &  7 & 10 & 12 & 13 & 17 & 19 & 20 & 23\\
\mathbf{15} &  2 &  3 &  4 &  5 &  8 & 11 & 12 & 16 & 18 & 21 & 24\\
\mathbf{16} &  1 &  5 &  6 &  8 &  9 & 11 & 12 & 15 & 18 & 21 & 25\\
\mathbf{17} &  2 &  6 &  7 &  8 & 10 & 12 & 13 & 14 & 19 & 20 & 26\\
\mathbf{18} &  2 &  4 &  8 &  9 & 10 & 11 & 12 & 15 & 16 & 21 & 27\\
\mathbf{19} &  1 &  3 &  5 &  6 &  7 & 10 & 13 & 14 & 17 & 20 & 28\\
\mathbf{20} &  2 &  3 &  6 &  7 &  8 & 10 & 13 & 14 & 17 & 19 & 29\\
\mathbf{21} &  1 &  4 &  5 &  8 &  9 & 11 & 12 & 15 & 16 & 18 & 30\\
\bottomrule
\end{array}\]

\subsection*{\rm Regular triangle-distinct graph \#$2$: $n=21$, $r=10$}

\[\begin{array}{c<{\hspace{1.5em}}cccccccccc>{\hspace{1.25em}}c}
\toprule
v & \multicolumn{10}{c}{\text{Neighbours of }v} & t(v)\\
\midrule
\mathbf{ 1} &  2 &  3 &  5 &  8 & 10 & 11 & 16 & 17 & 18 & 21 & 12\\
\mathbf{ 2} &  1 &  4 &  6 &  7 &  9 & 12 & 14 & 15 & 19 & 20 & 13\\
\mathbf{ 3} &  1 &  4 &  6 &  7 &  9 & 12 & 13 & 19 & 20 & 21 & 14\\
\mathbf{ 4} &  2 &  3 &  5 &  8 & 10 & 11 & 15 & 16 & 20 & 21 & 15\\
\mathbf{ 5} &  1 &  4 &  6 &  7 &  9 & 12 & 15 & 17 & 18 & 19 & 16\\
\mathbf{ 6} &  2 &  3 &  5 &  8 & 10 & 11 & 14 & 17 & 18 & 19 & 17\\
\mathbf{ 7} &  2 &  3 &  5 &  8 & 10 & 11 & 13 & 16 & 20 & 21 & 18\\
\mathbf{ 8} &  1 &  4 &  6 &  7 &  9 & 12 & 13 & 16 & 20 & 21 & 19\\
\mathbf{ 9} &  2 &  3 &  5 &  8 & 10 & 13 & 14 & 16 & 20 & 21 & 20\\
\mathbf{10} &  1 &  4 &  6 &  7 &  9 & 13 & 16 & 18 & 20 & 21 & 21\\
\mathbf{11} &  1 &  4 &  6 &  7 & 12 & 14 & 15 & 17 & 18 & 19 & 22\\
\mathbf{12} &  2 &  3 &  5 &  8 & 11 & 14 & 15 & 17 & 18 & 19 & 23\\
\mathbf{13} &  3 &  7 &  8 &  9 & 10 & 15 & 16 & 17 & 20 & 21 & 24\\
\mathbf{14} &  2 &  6 &  9 & 11 & 12 & 15 & 16 & 17 & 18 & 19 & 25\\
\mathbf{15} &  2 &  4 &  5 & 11 & 12 & 13 & 14 & 17 & 18 & 19 & 26\\
\mathbf{16} &  1 &  4 &  7 &  8 &  9 & 10 & 13 & 14 & 20 & 21 & 27\\
\mathbf{17} &  1 &  5 &  6 & 11 & 12 & 13 & 14 & 15 & 18 & 19 & 28\\
\mathbf{18} &  1 &  5 &  6 & 10 & 11 & 12 & 14 & 15 & 17 & 19 & 29\\
\mathbf{19} &  2 &  3 &  5 &  6 & 11 & 12 & 14 & 15 & 17 & 18 & 30\\
\mathbf{20} &  2 &  3 &  4 &  7 &  8 &  9 & 10 & 13 & 16 & 21 & 31\\
\mathbf{21} &  1 &  3 &  4 &  7 &  8 &  9 & 10 & 13 & 16 & 20 & 32\\
\bottomrule
\end{array}\]

\subsection*{\rm Regular triangle-distinct graph \#$3$: $n=21$, $r=10$}

\[\begin{array}{c<{\hspace{1.5em}}cccccccccc>{\hspace{1.25em}}c}
\toprule
v & \multicolumn{10}{c}{\text{Neighbours of }v} & t(v)\\
\midrule
\mathbf{ 1} &  2 &  3 &  4 &  7 &  9 & 12 & 13 & 14 & 20 & 21 & 11\\
\mathbf{ 2} &  1 &  5 &  6 &  8 & 10 & 11 & 16 & 17 & 18 & 21 & 12\\
\mathbf{ 3} &  1 &  5 &  6 &  8 & 10 & 11 & 15 & 17 & 18 & 21 & 13\\
\mathbf{ 4} &  1 &  5 &  6 &  8 & 10 & 11 & 16 & 18 & 19 & 20 & 14\\
\mathbf{ 5} &  2 &  3 &  4 &  7 &  9 & 12 & 14 & 15 & 19 & 20 & 15\\
\mathbf{ 6} &  2 &  3 &  4 &  9 & 12 & 13 & 14 & 16 & 19 & 20 & 16\\
\mathbf{ 7} &  1 &  5 &  8 & 11 & 13 & 15 & 16 & 17 & 18 & 21 & 17\\
\mathbf{ 8} &  2 &  3 &  4 &  7 & 13 & 14 & 15 & 16 & 19 & 20 & 18\\
\mathbf{ 9} &  1 &  5 &  6 & 10 & 11 & 13 & 17 & 18 & 19 & 21 & 19\\
\mathbf{10} &  2 &  3 &  4 &  9 & 12 & 13 & 15 & 17 & 18 & 21 & 20\\
\mathbf{11} &  2 &  3 &  4 &  7 &  9 & 13 & 14 & 17 & 18 & 21 & 21\\
\mathbf{12} &  1 &  5 &  6 & 10 & 14 & 15 & 16 & 17 & 19 & 20 & 22\\
\mathbf{13} &  1 &  6 &  7 &  8 &  9 & 10 & 11 & 17 & 18 & 21 & 23\\
\mathbf{14} &  1 &  5 &  6 &  8 & 11 & 12 & 15 & 16 & 19 & 20 & 24\\
\mathbf{15} &  3 &  5 &  7 &  8 & 10 & 12 & 14 & 16 & 19 & 20 & 25\\
\mathbf{16} &  2 &  4 &  6 &  7 &  8 & 12 & 14 & 15 & 19 & 20 & 26\\
\mathbf{17} &  2 &  3 &  7 &  9 & 10 & 11 & 12 & 13 & 18 & 21 & 27\\
\mathbf{18} &  2 &  3 &  4 &  7 &  9 & 10 & 11 & 13 & 17 & 21 & 28\\
\mathbf{19} &  4 &  5 &  6 &  8 &  9 & 12 & 14 & 15 & 16 & 20 & 29\\
\mathbf{20} &  1 &  4 &  5 &  6 &  8 & 12 & 14 & 15 & 16 & 19 & 30\\
\mathbf{21} &  1 &  2 &  3 &  7 &  9 & 10 & 11 & 13 & 17 & 18 & 31\\
\bottomrule
\end{array}\]

\subsection*{\rm Regular triangle-distinct graph \#$4$: $n=21$, $r=10$}

\[\begin{array}{c<{\hspace{1.5em}}cccccccccc>{\hspace{1.25em}}c}
\toprule
v & \multicolumn{10}{c}{\text{Neighbours of }v} & t(v)\\
\midrule
\mathbf{ 1} &  2 &  3 &  6 &  8 &  9 & 11 & 12 & 13 & 18 & 21 &  9\\
\mathbf{ 2} &  1 &  4 &  5 &  7 & 10 & 14 & 15 & 16 & 19 & 20 & 10\\
\mathbf{ 3} &  1 &  4 &  5 &  7 & 10 & 14 & 15 & 16 & 20 & 21 & 11\\
\mathbf{ 4} &  2 &  3 &  6 &  9 & 11 & 12 & 13 & 17 & 18 & 21 & 13\\
\mathbf{ 5} &  2 &  3 &  6 &  8 &  9 & 11 & 12 & 18 & 19 & 20 & 14\\
\mathbf{ 6} &  1 &  4 &  5 &  7 & 10 & 13 & 14 & 15 & 19 & 20 & 15\\
\mathbf{ 7} &  2 &  3 &  6 & 11 & 12 & 13 & 16 & 17 & 18 & 21 & 16\\
\mathbf{ 8} &  1 &  5 & 10 & 12 & 14 & 15 & 16 & 17 & 19 & 20 & 17\\
\mathbf{ 9} &  1 &  4 &  5 & 10 & 15 & 16 & 17 & 18 & 19 & 21 & 18\\
\mathbf{10} &  2 &  3 &  6 &  8 &  9 & 11 & 16 & 17 & 18 & 21 & 19\\
\mathbf{11} &  1 &  4 &  5 &  7 & 10 & 14 & 16 & 17 & 18 & 21 & 20\\
\mathbf{12} &  1 &  4 &  5 &  7 &  8 & 13 & 14 & 15 & 19 & 20 & 21\\
\mathbf{13} &  1 &  4 &  6 &  7 & 12 & 14 & 15 & 17 & 19 & 20 & 22\\
\mathbf{14} &  2 &  3 &  6 &  8 & 11 & 12 & 13 & 15 & 19 & 20 & 23\\
\mathbf{15} &  2 &  3 &  6 &  8 &  9 & 12 & 13 & 14 & 19 & 20 & 24\\
\mathbf{16} &  2 &  3 &  7 &  8 &  9 & 10 & 11 & 17 & 18 & 21 & 25\\
\mathbf{17} &  4 &  7 &  8 &  9 & 10 & 11 & 13 & 16 & 18 & 21 & 26\\
\mathbf{18} &  1 &  4 &  5 &  7 &  9 & 10 & 11 & 16 & 17 & 21 & 27\\
\mathbf{19} &  2 &  5 &  6 &  8 &  9 & 12 & 13 & 14 & 15 & 20 & 28\\
\mathbf{20} &  2 &  3 &  5 &  6 &  8 & 12 & 13 & 14 & 15 & 19 & 29\\
\mathbf{21} &  1 &  3 &  4 &  7 &  9 & 10 & 11 & 16 & 17 & 18 & 30\\
\bottomrule
\end{array}\]

\subsection*{\rm Regular triangle-distinct graph \#$5$: $n=22$, $r=10$}

\[\begin{array}{c<{\hspace{1.5em}}cccccccccc>{\hspace{1.25em}}c}
\toprule
v & \multicolumn{10}{c}{\text{Neighbours of }v} & t(v)\\
\midrule
\mathbf{ 1} &  4 &  5 &  6 &  7 &  9 & 10 & 11 & 14 & 18 & 22 &  9\\
\mathbf{ 2} &  4 &  5 &  6 &  7 &  9 & 10 & 11 & 14 & 18 & 21 & 10\\
\mathbf{ 3} &  4 &  5 &  6 &  7 &  9 & 10 & 11 & 14 & 16 & 19 & 11\\
\mathbf{ 4} &  1 &  2 &  3 &  8 & 12 & 15 & 16 & 17 & 19 & 21 & 12\\
\mathbf{ 5} &  1 &  2 &  3 &  8 & 12 & 16 & 17 & 19 & 20 & 21 & 13\\
\mathbf{ 6} &  1 &  2 &  3 &  8 & 12 & 13 & 15 & 16 & 20 & 22 & 14\\
\mathbf{ 7} &  1 &  2 &  3 &  8 & 12 & 13 & 15 & 18 & 20 & 22 & 15\\
\mathbf{ 8} &  4 &  5 &  6 &  7 &  9 & 10 & 17 & 20 & 21 & 22 & 16\\
\mathbf{ 9} &  1 &  2 &  3 &  8 & 15 & 16 & 17 & 18 & 20 & 22 & 17\\
\mathbf{10} &  1 &  2 &  3 &  8 & 12 & 13 & 14 & 17 & 19 & 21 & 18\\
\mathbf{11} &  1 &  2 &  3 & 12 & 13 & 14 & 17 & 18 & 19 & 21 & 19\\
\mathbf{12} &  4 &  5 &  6 &  7 & 10 & 11 & 14 & 17 & 19 & 21 & 20\\
\mathbf{13} &  6 &  7 & 10 & 11 & 15 & 16 & 18 & 19 & 20 & 22 & 21\\
\mathbf{14} &  1 &  2 &  3 & 10 & 11 & 12 & 15 & 17 & 19 & 21 & 22\\
\mathbf{15} &  4 &  6 &  7 &  9 & 13 & 14 & 16 & 18 & 20 & 22 & 23\\
\mathbf{16} &  3 &  4 &  5 &  6 &  9 & 13 & 15 & 18 & 20 & 22 & 24\\
\mathbf{17} &  4 &  5 &  8 &  9 & 10 & 11 & 12 & 14 & 19 & 21 & 25\\
\mathbf{18} &  1 &  2 &  7 &  9 & 11 & 13 & 15 & 16 & 20 & 22 & 26\\
\mathbf{19} &  3 &  4 &  5 & 10 & 11 & 12 & 13 & 14 & 17 & 21 & 27\\
\mathbf{20} &  5 &  6 &  7 &  8 &  9 & 13 & 15 & 16 & 18 & 22 & 28\\
\mathbf{21} &  2 &  4 &  5 &  8 & 10 & 11 & 12 & 14 & 17 & 19 & 29\\
\mathbf{22} &  1 &  6 &  7 &  8 &  9 & 13 & 15 & 16 & 18 & 20 & 30\\
\bottomrule
\end{array}\]

\subsection*{\rm Regular triangle-distinct graph \#$6$: $n=23$, $r=10$}

\[\begin{array}{c<{\hspace{1.5em}}cccccccccc>{\hspace{1.25em}}c}
\toprule
v & \multicolumn{10}{c}{\text{Neighbours of }v} & t(v)\\
\midrule
\mathbf{ 1} &  3 &  4 &  5 &  9 & 10 & 11 & 12 & 16 & 18 & 23 &  7\\
\mathbf{ 2} &  3 &  4 &  5 &  9 & 10 & 11 & 12 & 18 & 19 & 21 &  9\\
\mathbf{ 3} &  1 &  2 &  6 &  7 &  8 & 13 & 14 & 15 & 21 & 22 & 10\\
\mathbf{ 4} &  1 &  2 &  6 &  7 &  8 & 13 & 15 & 20 & 21 & 22 & 11\\
\mathbf{ 5} &  1 &  2 &  6 &  8 & 13 & 14 & 15 & 19 & 20 & 23 & 12\\
\mathbf{ 6} &  3 &  4 &  5 &  9 & 10 & 11 & 12 & 18 & 21 & 22 & 13\\
\mathbf{ 7} &  3 &  4 &  9 & 10 & 11 & 12 & 13 & 16 & 18 & 22 & 14\\
\mathbf{ 8} &  3 &  4 &  5 & 11 & 12 & 16 & 18 & 19 & 21 & 22 & 15\\
\mathbf{ 9} &  1 &  2 &  6 &  7 & 14 & 15 & 17 & 19 & 20 & 23 & 16\\
\mathbf{10} &  1 &  2 &  6 &  7 & 14 & 15 & 16 & 17 & 20 & 23 & 17\\
\mathbf{11} &  1 &  2 &  6 &  7 &  8 & 13 & 14 & 18 & 21 & 22 & 18\\
\mathbf{12} &  1 &  2 &  6 &  7 &  8 & 13 & 17 & 18 & 21 & 22 & 19\\
\mathbf{13} &  3 &  4 &  5 &  7 & 11 & 12 & 17 & 18 & 21 & 22 & 20\\
\mathbf{14} &  3 &  5 &  9 & 10 & 11 & 16 & 17 & 19 & 20 & 23 & 21\\
\mathbf{15} &  3 &  4 &  5 &  9 & 10 & 16 & 17 & 19 & 20 & 23 & 22\\
\mathbf{16} &  1 &  7 &  8 & 10 & 14 & 15 & 17 & 19 & 20 & 23 & 23\\
\mathbf{17} &  9 & 10 & 12 & 13 & 14 & 15 & 16 & 19 & 20 & 23 & 25\\
\mathbf{18} &  1 &  2 &  6 &  7 &  8 & 11 & 12 & 13 & 21 & 22 & 26\\
\mathbf{19} &  2 &  5 &  8 &  9 & 14 & 15 & 16 & 17 & 20 & 23 & 27\\
\mathbf{20} &  4 &  5 &  9 & 10 & 14 & 15 & 16 & 17 & 19 & 23 & 29\\
\mathbf{21} &  2 &  3 &  4 &  6 &  8 & 11 & 12 & 13 & 18 & 22 & 30\\
\mathbf{22} &  3 &  4 &  6 &  7 &  8 & 11 & 12 & 13 & 18 & 21 & 31\\
\mathbf{23} &  1 &  5 &  9 & 10 & 14 & 15 & 16 & 17 & 19 & 20 & 32\\
\bottomrule
\end{array}\]

\subsection*{\rm Regular triangle-distinct graph \#$7$: $n=24$, $r=11$}

\[\begin{array}{c<{\hspace{1.5em}}ccccccccccc>{\hspace{1.25em}}c}
\toprule
v & \multicolumn{11}{c}{\text{Neighbours of }v} & t(v)\\
\midrule
\mathbf{ 1} &  2 &  3 &  8 &  9 & 11 & 12 & 13 & 15 & 17 & 20 & 21 & 13\\
\mathbf{ 2} &  1 &  4 &  5 &  6 &  7 & 10 & 16 & 18 & 19 & 21 & 22 & 14\\
\mathbf{ 3} &  1 &  4 &  5 &  6 &  7 & 10 & 16 & 18 & 22 & 23 & 24 & 15\\
\mathbf{ 4} &  2 &  3 &  8 &  9 & 10 & 11 & 12 & 13 & 20 & 22 & 23 & 16\\
\mathbf{ 5} &  2 &  3 &  9 & 10 & 11 & 12 & 13 & 17 & 20 & 21 & 24 & 17\\
\mathbf{ 6} &  2 &  3 &  8 &  9 & 11 & 13 & 14 & 15 & 20 & 23 & 24 & 18\\
\mathbf{ 7} &  2 &  3 &  8 &  9 & 12 & 13 & 15 & 20 & 21 & 23 & 24 & 19\\
\mathbf{ 8} &  1 &  4 &  6 &  7 & 10 & 15 & 16 & 17 & 18 & 19 & 22 & 20\\
\mathbf{ 9} &  1 &  4 &  5 &  6 &  7 & 12 & 16 & 17 & 18 & 19 & 22 & 21\\
\mathbf{10} &  2 &  3 &  4 &  5 &  8 & 14 & 17 & 19 & 21 & 23 & 24 & 22\\
\mathbf{11} &  1 &  4 &  5 &  6 & 12 & 14 & 16 & 17 & 18 & 19 & 22 & 23\\
\mathbf{12} &  1 &  4 &  5 &  7 &  9 & 11 & 14 & 16 & 18 & 19 & 22 & 24\\
\mathbf{13} &  1 &  4 &  5 &  6 &  7 & 14 & 16 & 20 & 21 & 23 & 24 & 25\\
\mathbf{14} &  6 & 10 & 11 & 12 & 13 & 15 & 18 & 20 & 21 & 23 & 24 & 26\\
\mathbf{15} &  1 &  6 &  7 &  8 & 14 & 17 & 19 & 20 & 21 & 23 & 24 & 27\\
\mathbf{16} &  2 &  3 &  8 &  9 & 11 & 12 & 13 & 17 & 18 & 19 & 22 & 28\\
\mathbf{17} &  1 &  5 &  8 &  9 & 10 & 11 & 15 & 16 & 18 & 19 & 22 & 29\\
\mathbf{18} &  2 &  3 &  8 &  9 & 11 & 12 & 14 & 16 & 17 & 19 & 22 & 30\\
\mathbf{19} &  2 &  8 &  9 & 10 & 11 & 12 & 15 & 16 & 17 & 18 & 22 & 31\\
\mathbf{20} &  1 &  4 &  5 &  6 &  7 & 13 & 14 & 15 & 21 & 23 & 24 & 32\\
\mathbf{21} &  1 &  2 &  5 &  7 & 10 & 13 & 14 & 15 & 20 & 23 & 24 & 33\\
\mathbf{22} &  2 &  3 &  4 &  8 &  9 & 11 & 12 & 16 & 17 & 18 & 19 & 34\\
\mathbf{23} &  3 &  4 &  6 &  7 & 10 & 13 & 14 & 15 & 20 & 21 & 24 & 35\\
\mathbf{24} &  3 &  5 &  6 &  7 & 10 & 13 & 14 & 15 & 20 & 21 & 23 & 36\\
\bottomrule
\end{array}\]

\subsection*{\rm Regular triangle-distinct graph \#$8$: $n=24$, $r=11$}

\[\begin{array}{c<{\hspace{1.5em}}ccccccccccc>{\hspace{1.25em}}c}
\toprule
v & \multicolumn{11}{c}{\text{Neighbours of }v} & t(v)\\
\midrule
\mathbf{ 1} &  2 &  4 &  5 &  9 & 11 & 12 & 14 & 17 & 18 & 21 & 23 & 12\\
\mathbf{ 2} &  1 &  3 &  6 &  7 &  8 & 10 & 13 & 15 & 16 & 20 & 23 & 13\\
\mathbf{ 3} &  2 &  4 &  5 &  9 & 11 & 12 & 17 & 18 & 21 & 23 & 24 & 14\\
\mathbf{ 4} &  1 &  3 &  6 &  7 &  8 & 10 & 15 & 16 & 19 & 22 & 24 & 15\\
\mathbf{ 5} &  1 &  3 &  6 &  7 &  8 & 10 & 13 & 14 & 16 & 22 & 24 & 16\\
\mathbf{ 6} &  2 &  4 &  5 &  9 & 12 & 14 & 17 & 18 & 19 & 21 & 23 & 17\\
\mathbf{ 7} &  2 &  4 &  5 &  9 & 11 & 12 & 14 & 15 & 18 & 19 & 20 & 18\\
\mathbf{ 8} &  2 &  4 &  5 &  9 & 11 & 12 & 13 & 14 & 17 & 22 & 24 & 19\\
\mathbf{ 9} &  1 &  3 &  6 &  7 &  8 & 13 & 15 & 16 & 17 & 22 & 24 & 20\\
\mathbf{10} &  2 &  4 &  5 & 11 & 16 & 17 & 18 & 20 & 21 & 22 & 24 & 21\\
\mathbf{11} &  1 &  3 &  7 &  8 & 10 & 15 & 16 & 17 & 19 & 22 & 24 & 22\\
\mathbf{12} &  1 &  3 &  6 &  7 &  8 & 16 & 18 & 19 & 20 & 21 & 23 & 23\\
\mathbf{13} &  2 &  5 &  8 &  9 & 14 & 18 & 19 & 20 & 21 & 22 & 23 & 24\\
\mathbf{14} &  1 &  5 &  6 &  7 &  8 & 13 & 15 & 19 & 20 & 21 & 23 & 25\\
\mathbf{15} &  2 &  4 &  7 &  9 & 11 & 14 & 16 & 17 & 20 & 22 & 24 & 26\\
\mathbf{16} &  2 &  4 &  5 &  9 & 10 & 11 & 12 & 15 & 17 & 22 & 24 & 27\\
\mathbf{17} &  1 &  3 &  6 &  8 &  9 & 10 & 11 & 15 & 16 & 22 & 24 & 28\\
\mathbf{18} &  1 &  3 &  6 &  7 & 10 & 12 & 13 & 19 & 20 & 21 & 23 & 29\\
\mathbf{19} &  4 &  6 &  7 & 11 & 12 & 13 & 14 & 18 & 20 & 21 & 23 & 30\\
\mathbf{20} &  2 &  7 & 10 & 12 & 13 & 14 & 15 & 18 & 19 & 21 & 23 & 31\\
\mathbf{21} &  1 &  3 &  6 & 10 & 12 & 13 & 14 & 18 & 19 & 20 & 23 & 32\\
\mathbf{22} &  4 &  5 &  8 &  9 & 10 & 11 & 13 & 15 & 16 & 17 & 24 & 34\\
\mathbf{23} &  1 &  2 &  3 &  6 & 12 & 13 & 14 & 18 & 19 & 20 & 21 & 35\\
\mathbf{24} &  3 &  4 &  5 &  8 &  9 & 10 & 11 & 15 & 16 & 17 & 22 & 36\\
\bottomrule
\end{array}\]

\subsection*{\rm Regular triangle-distinct graph \#$9$: $n=25$, $r=12$}

\[\begin{array}{c<{\hspace{1.5em}}cccccccccccc>{\hspace{1.25em}}c}
\toprule
v & \multicolumn{12}{c}{\text{Neighbours of }v} & t(v)\\
\midrule
\mathbf{ 1} &  2 &  3 &  4 & 10 & 11 & 12 & 14 & 16 & 17 & 18 & 21 & 24 & 18\\
\mathbf{ 2} &  1 &  5 &  6 &  7 &  8 &  9 & 13 & 15 & 20 & 22 & 23 & 24 & 19\\
\mathbf{ 3} &  1 &  5 &  6 &  7 &  8 &  9 & 13 & 17 & 19 & 20 & 22 & 25 & 20\\
\mathbf{ 4} &  1 &  5 &  6 &  7 &  8 &  9 & 13 & 17 & 19 & 22 & 24 & 25 & 21\\
\mathbf{ 5} &  2 &  3 &  4 & 10 & 11 & 12 & 14 & 16 & 17 & 19 & 23 & 25 & 22\\
\mathbf{ 6} &  2 &  3 &  4 & 10 & 11 & 12 & 14 & 16 & 18 & 19 & 22 & 25 & 23\\
\mathbf{ 7} &  2 &  3 &  4 & 10 & 11 & 14 & 16 & 18 & 20 & 21 & 23 & 24 & 24\\
\mathbf{ 8} &  2 &  3 &  4 & 10 & 11 & 12 & 14 & 16 & 17 & 18 & 22 & 25 & 25\\
\mathbf{ 9} &  2 &  3 &  4 & 11 & 12 & 14 & 15 & 18 & 19 & 21 & 23 & 24 & 26\\
\mathbf{10} &  1 &  5 &  6 &  7 &  8 & 11 & 13 & 15 & 16 & 20 & 22 & 25 & 27\\
\mathbf{11} &  1 &  5 &  6 &  7 &  8 &  9 & 10 & 13 & 14 & 15 & 22 & 25 & 28\\
\mathbf{12} &  1 &  5 &  6 &  8 &  9 & 13 & 14 & 17 & 21 & 22 & 23 & 25 & 29\\
\mathbf{13} &  2 &  3 &  4 & 10 & 11 & 12 & 14 & 15 & 17 & 20 & 22 & 25 & 31\\
\mathbf{14} &  1 &  5 &  6 &  7 &  8 &  9 & 11 & 12 & 13 & 17 & 22 & 25 & 32\\
\mathbf{15} &  2 &  9 & 10 & 11 & 13 & 17 & 18 & 19 & 20 & 21 & 23 & 24 & 33\\
\mathbf{16} &  1 &  5 &  6 &  7 &  8 & 10 & 18 & 19 & 20 & 21 & 23 & 24 & 34\\
\mathbf{17} &  1 &  3 &  4 &  5 &  8 & 12 & 13 & 14 & 15 & 21 & 22 & 25 & 35\\
\mathbf{18} &  1 &  6 &  7 &  8 &  9 & 15 & 16 & 19 & 20 & 21 & 23 & 24 & 36\\
\mathbf{19} &  3 &  4 &  5 &  6 &  9 & 15 & 16 & 18 & 20 & 21 & 23 & 24 & 37\\
\mathbf{20} &  2 &  3 &  7 & 10 & 13 & 15 & 16 & 18 & 19 & 21 & 23 & 24 & 38\\
\mathbf{21} &  1 &  7 &  9 & 12 & 15 & 16 & 17 & 18 & 19 & 20 & 23 & 24 & 39\\
\mathbf{22} &  2 &  3 &  4 &  6 &  8 & 10 & 11 & 12 & 13 & 14 & 17 & 25 & 40\\
\mathbf{23} &  2 &  5 &  7 &  9 & 12 & 15 & 16 & 18 & 19 & 20 & 21 & 24 & 41\\
\mathbf{24} &  1 &  2 &  4 &  7 &  9 & 15 & 16 & 18 & 19 & 20 & 21 & 23 & 43\\
\mathbf{25} &  3 &  4 &  5 &  6 &  8 & 10 & 11 & 12 & 13 & 14 & 17 & 22 & 44\\
\bottomrule
\end{array}\]

\subsection*{\rm Regular triangle-distinct graph \#$10$: $n=26$, $r=12$}

\[\begin{array}{c<{\hspace{1.5em}}cccccccccccc>{\hspace{1.25em}}c}
\toprule
v & \multicolumn{12}{c}{\text{Neighbours of }v} & t(v)\\
\midrule
\mathbf{ 1} &  2 &  3 &  4 &  7 &  9 & 12 & 13 & 14 & 15 & 19 & 23 & 25 & 13\\
\mathbf{ 2} &  1 &  5 &  6 &  8 & 10 & 11 & 16 & 17 & 18 & 20 & 21 & 22 & 14\\
\mathbf{ 3} &  1 &  5 &  6 &  8 & 10 & 11 & 17 & 18 & 20 & 21 & 24 & 26 & 16\\
\mathbf{ 4} &  1 &  5 &  6 &  8 & 10 & 11 & 16 & 17 & 18 & 20 & 22 & 26 & 17\\
\mathbf{ 5} &  2 &  3 &  4 &  7 &  9 & 12 & 13 & 14 & 19 & 23 & 24 & 25 & 18\\
\mathbf{ 6} &  2 &  3 &  4 &  7 &  9 & 12 & 13 & 14 & 19 & 21 & 24 & 25 & 19\\
\mathbf{ 7} &  1 &  5 &  6 &  8 & 10 & 15 & 16 & 17 & 21 & 22 & 24 & 26 & 20\\
\mathbf{ 8} &  2 &  3 &  4 &  7 &  9 & 12 & 14 & 15 & 18 & 23 & 25 & 26 & 21\\
\mathbf{ 9} &  1 &  5 &  6 &  8 & 11 & 18 & 20 & 22 & 23 & 24 & 25 & 26 & 22\\
\mathbf{10} &  2 &  3 &  4 &  7 & 12 & 13 & 15 & 19 & 20 & 21 & 22 & 23 & 23\\
\mathbf{11} &  2 &  3 &  4 &  9 & 13 & 14 & 15 & 16 & 17 & 19 & 24 & 26 & 24\\
\mathbf{12} &  1 &  5 &  6 &  8 & 10 & 16 & 18 & 20 & 21 & 23 & 24 & 25 & 25\\
\mathbf{13} &  1 &  5 &  6 & 10 & 11 & 14 & 17 & 18 & 19 & 22 & 23 & 25 & 26\\
\mathbf{14} &  1 &  5 &  6 &  8 & 11 & 13 & 15 & 17 & 18 & 20 & 22 & 26 & 27\\
\mathbf{15} &  1 &  7 &  8 & 10 & 11 & 14 & 16 & 18 & 20 & 21 & 22 & 26 & 28\\
\mathbf{16} &  2 &  4 &  7 & 11 & 12 & 15 & 17 & 19 & 21 & 23 & 24 & 25 & 29\\
\mathbf{17} &  2 &  3 &  4 &  7 & 11 & 13 & 14 & 16 & 19 & 20 & 22 & 26 & 30\\
\mathbf{18} &  2 &  3 &  4 &  8 &  9 & 12 & 13 & 14 & 15 & 20 & 22 & 26 & 31\\
\mathbf{19} &  1 &  5 &  6 & 10 & 11 & 13 & 16 & 17 & 21 & 23 & 24 & 25 & 32\\
\mathbf{20} &  2 &  3 &  4 &  9 & 10 & 12 & 14 & 15 & 17 & 18 & 22 & 26 & 33\\
\mathbf{21} &  2 &  3 &  6 &  7 & 10 & 12 & 15 & 16 & 19 & 23 & 24 & 25 & 34\\
\mathbf{22} &  2 &  4 &  7 &  9 & 10 & 13 & 14 & 15 & 17 & 18 & 20 & 26 & 35\\
\mathbf{23} &  1 &  5 &  8 &  9 & 10 & 12 & 13 & 16 & 19 & 21 & 24 & 25 & 36\\
\mathbf{24} &  3 &  5 &  6 &  7 &  9 & 11 & 12 & 16 & 19 & 21 & 23 & 25 & 37\\
\mathbf{25} &  1 &  5 &  6 &  8 &  9 & 12 & 13 & 16 & 19 & 21 & 23 & 24 & 38\\
\mathbf{26} &  3 &  4 &  7 &  8 &  9 & 11 & 14 & 15 & 17 & 18 & 20 & 22 & 39\\
\bottomrule
\end{array}\]

\subsection*{\rm Regular triangle-distinct graph \#$11$: $n=27$, $r=12$}

\[\begin{array}{c<{\hspace{1.5em}}cccccccccccc>{\hspace{1.25em}}c}
\toprule
v & \multicolumn{12}{c}{\text{Neighbours of }v} & t(v)\\
\midrule
\mathbf{ 1} &  3 &  4 &  5 &  6 &  7 & 12 & 14 & 15 & 17 & 20 & 21 & 24 & 14\\
\mathbf{ 2} &  3 &  4 &  5 &  6 &  7 & 12 & 14 & 15 & 17 & 21 & 24 & 25 & 15\\
\mathbf{ 3} &  1 &  2 &  8 &  9 & 10 & 11 & 13 & 16 & 18 & 19 & 26 & 27 & 17\\
\mathbf{ 4} &  1 &  2 &  8 & 10 & 13 & 16 & 18 & 19 & 21 & 22 & 23 & 27 & 19\\
\mathbf{ 5} &  1 &  2 &  9 & 10 & 11 & 16 & 18 & 19 & 20 & 22 & 23 & 26 & 20\\
\mathbf{ 6} &  1 &  2 &  8 &  9 & 10 & 13 & 17 & 18 & 19 & 22 & 26 & 27 & 21\\
\mathbf{ 7} &  1 &  2 &  9 & 10 & 11 & 13 & 16 & 19 & 20 & 23 & 24 & 25 & 22\\
\mathbf{ 8} &  3 &  4 &  6 & 12 & 14 & 15 & 18 & 20 & 24 & 25 & 26 & 27 & 24\\
\mathbf{ 9} &  3 &  5 &  6 &  7 & 12 & 14 & 15 & 21 & 22 & 23 & 24 & 25 & 25\\
\mathbf{10} &  3 &  4 &  5 &  6 &  7 & 12 & 15 & 18 & 22 & 23 & 24 & 25 & 26\\
\mathbf{11} &  3 &  5 &  7 & 12 & 15 & 16 & 17 & 20 & 21 & 22 & 23 & 25 & 27\\
\mathbf{12} &  1 &  2 &  8 &  9 & 10 & 11 & 17 & 18 & 21 & 22 & 24 & 25 & 28\\
\mathbf{13} &  3 &  4 &  6 &  7 & 14 & 15 & 16 & 17 & 19 & 21 & 26 & 27 & 29\\
\mathbf{14} &  1 &  2 &  8 &  9 & 13 & 16 & 17 & 18 & 19 & 20 & 26 & 27 & 30\\
\mathbf{15} &  1 &  2 &  8 &  9 & 10 & 11 & 13 & 21 & 22 & 23 & 24 & 25 & 31\\
\mathbf{16} &  3 &  4 &  5 &  7 & 11 & 13 & 14 & 17 & 20 & 23 & 26 & 27 & 32\\
\mathbf{17} &  1 &  2 &  6 & 11 & 12 & 13 & 14 & 16 & 19 & 20 & 26 & 27 & 33\\
\mathbf{18} &  3 &  4 &  5 &  6 &  8 & 10 & 12 & 14 & 19 & 20 & 26 & 27 & 34\\
\mathbf{19} &  3 &  4 &  5 &  6 &  7 & 13 & 14 & 17 & 18 & 20 & 26 & 27 & 35\\
\mathbf{20} &  1 &  5 &  7 &  8 & 11 & 14 & 16 & 17 & 18 & 19 & 26 & 27 & 36\\
\mathbf{21} &  1 &  2 &  4 &  9 & 11 & 12 & 13 & 15 & 22 & 23 & 24 & 25 & 37\\
\mathbf{22} &  4 &  5 &  6 &  9 & 10 & 11 & 12 & 15 & 21 & 23 & 24 & 25 & 38\\
\mathbf{23} &  4 &  5 &  7 &  9 & 10 & 11 & 15 & 16 & 21 & 22 & 24 & 25 & 39\\
\mathbf{24} &  1 &  2 &  7 &  8 &  9 & 10 & 12 & 15 & 21 & 22 & 23 & 25 & 40\\
\mathbf{25} &  2 &  7 &  8 &  9 & 10 & 11 & 12 & 15 & 21 & 22 & 23 & 24 & 42\\
\mathbf{26} &  3 &  5 &  6 &  8 & 13 & 14 & 16 & 17 & 18 & 19 & 20 & 27 & 43\\
\mathbf{27} &  3 &  4 &  6 &  8 & 13 & 14 & 16 & 17 & 18 & 19 & 20 & 26 & 44\\
\bottomrule
\end{array}\]

\end{document}